\documentclass[letterpaper,10pt]{amsart}

\usepackage[latin1]{inputenc}
\usepackage{amsfonts}
\usepackage{amsmath}
\usepackage{amssymb}
\usepackage{amsthm}
\usepackage{amsrefs}
\usepackage{hyperref}
\usepackage{empheq}
\usepackage{nicefrac}
\usepackage{enumitem}  
\usepackage[all]{xy}
\usepackage{graphicx}
\usepackage{textcomp}

\usepackage{mathtools}
\usepackage{bbm}

\theoremstyle{plain}

\theoremstyle{definition}
	\newtheorem*{theorem*}{Theorem}
	\newtheorem{theorem}{Theorem}
	\newtheorem*{prop*}{Proposition}
	\newtheorem{prop}{Proposition}
	\newtheorem{remark}[prop]{Remark}
	\newtheorem*{remark*}{Remark}	
	\newtheorem*{lemma*}{Lemma}
	\newtheorem{lemma}[prop]{Lemma}

	\newtheorem*{definition*}{Definition}
	
	\newtheorem*{question*}{Question}
	
	\newtheorem*{note*}{Note}
	
	\newtheorem*{claim*}{Claim}
	\newtheorem*{corollary*}{Corollary}
	
	\newtheorem*{example*}{Example}

\numberwithin{prop}{section}
\numberwithin{corollary}{section}
\numberwithin{example}{section}	
\numberwithin{equation}{subsection} 
\numberwithin{claim}{section}	

\newcommand{\Z}{\mathbb{Z}}
\newcommand{\R}{\mathbb{R}}

\newcommand{\C}{\mathbb{C}}

\begin{document}

\title[The Automorphic Heat Kernel]{Global Automorphic Sobolev Theory \\and The~Automorphic Heat~Kernel}

\author{Amy T. DeCelles}
\address{University of St. Thomas \\2115 Summit Avenue \\
St. Paul, MN 55105}
\email{adecelles@stthomas.edu}
\urladdr{http://personal.stthomas.edu/dece4515}


\thanks{The idea for this paper was suggested to the author by Paul Garrett.  The author would like to thank him and the other members of the Automorphic Forms Seminar at the University of Minnesota for helpful conversations.  The author was partially supported by a research grant from the University of St. Thomas.}

\begin{abstract} Heat kernels arise in a variety of contexts including probability, geometry, and functional analysis; the automorphic heat kernel is particularly important in number theory and string theory.  The typical construction of an automorphic heat kernel as a Poincar\'{e} series presents analytic difficulties, which can be dealt with in special cases (e.g. hyperbolic spaces) but are often sidestepped in higher rank by restricting to the compact quotient case.  In this paper, we present a new approach, using global automorphic Sobolev theory, a robust framework for solving automorphic PDEs that does not require any simplifying assumptions about the rank of the symmetric space or the compactness of the arithmetic quotient.  We construct an automorphic heat kernel via its automorphic spectral expansion in terms of cusp forms, Eisenstein series, and residues of Eisenstein series.  We then prove uniqueness of the automorphic heat kernel as an application of operator semigroup theory.  Finally, we prove the smoothness of the automorphic heat kernel by proving that its automorphic spectral expansion converges in the $C^\infty$-topology.  
\end{abstract}

\maketitle

\section{Introduction}
Heat kernels arise in many contexts in mathematics and physics, so many that Jorgenson and Lang have called the heat kernel ``ubiquitous'' \cite{jorgenson-lang2001}.  Automorphic heat kernels are important in physics and number theory: applications include  multiloop amplitudes of certain bosonic strings \cite{dhoker-phong1986}, asymptotic formulas for spectra of arithmetic quotients \cite{gangolli1968, donnelly1982, deitmar-hoffman1999}, a relationship between eta invariants of certain manifolds and closed geodesics using the Selberg trace formula \cite{moscovici-stanton1989}, systematic construction of zeta-type functions via heat Eisenstein series \cite{jorgenson-lang2001, jorgenson-lang2008, jorgenson-lang2009, jorgenson-lang2012}, sup-norm bounds for automorphic forms \cite{jorgenson-kramer2004, jorgenson-kramer2011, ary2016, friedman-jorgenson-kramer2016, ary-bal2018}, limit formulas for period integrals and a Weyl-type asymptotic law for a counting function for period integrals, via identities that can be considered as a special cases of Jacquet's relative trace formula  \cite{tsuzuki2008, tsuzuki2009}, and an average version of the holomorphic QUE conjecture for automorphic cusp forms associated to quaternion algebras \cite{ary-bal2018}.  Heat asymptotics on spaces of automorphic forms are of continuing interest; see \cite{paniagua-taboada-ramacher2016} and its references.

Typically, an automorphic heat kernel is constructed as a Poincar\'{e} series, obtained by ``winding up'' or ``periodicizing'' a heat kernel on a symmetric space $G/K$ by averaging over an arithmetic subgroup $\Gamma \subset G$.  Proving convergence of the resulting series is a non-trivial issue.  In special cases, e.g. when $G$ is of rank one or complex, convergence may be proven using explicit computations \cite{fay1977, oshima1990, jorgenson-lang2008, jorgenson-lang2009}.  Many of the applications listed above take advantage of such cases.  Other applications are obtained when $\Gamma \backslash G$ is compact; sometimes this assumption is explicit as in \cite{moscovici-stanton1989, tsuzuki2008, tsuzuki2009}; other times results are stated in generality but convergence proofs seem to assume $\Gamma \backslash G$ is compact, as in \cite{donnelly1982, deitmar-hoffman1999}, relying on \cite{donnelly1979}, which treats only the compact quotient case.  Another issue with the Poincar\'{e} series construction of an automorphic heat kernel is the difficulty of obtaining an automorphic spectral expansion.  In particular, this seems to be a significant obstacle in Jorgenson and Lang's program for systematic construction of zeta-type functions from heat Eisenstein series; see Sections 4.6 and 4.7 of \cite{jorgenson-lang2009}.

In this paper, we take an entirely different approach: we construct an automorphic heat kernel on $\Gamma \backslash G/K$ via its automorphic spectral expansion in terms of cusp forms, Eisenstein series, and residues of Eisenstein series.  Here $G$ is of arbitrary rank, and $\Gamma \backslash G$ may be noncompact.  Having proven existence of an automorphic heat kernel via its spectral expansion, we then prove uniqueness and smoothness.  Admittedly, the existence, uniqueness, and smoothness of the automorphic heat kernel are not surprising, and the automorphic spectral expansion is exactly what heuristic arguments predict.  What is novel is the framework in which these results are proven.  Moreover it is  noteworthy that, once in this framework, the results are proven decisively in such generality.  

This framework, developed in \cite{decelles2012}, is global automorphic Sobolev theory, which allows us to solve automorphic PDEs by translating to the spectral side and back again in a completely rigorous way.  In this paper, the main theorems about the automorphic heat kernel are all obtained as easy corollaries of the corresponding results on the simpler, spectral side, where the Laplacian is replaced by a multiplication operator (or a family of multiplication operators).  

For other applications of global automorphic Sobolev theory to number theory and physics, see \cite{decelles2012, decelles2016, klinger-logan2018}; applications include lattice point counting in symmetric spaces and the behavior the 4-loop supergraviton.  The presence of the time parameter in the heat equation adds a layer of complexity not found in these other applications; to construct the automorphic heat kernel, we solve a \emph{vector-valued} differential equation on the spectral side.  We need the theory of differentiation of vector-valued functions due to Groethendick and Schwartz as well as the theory of strongly continuous semigroups of linear operators.

\subsection*{Preliminaries.} Let $G$ be a reductive or semi-simple Lie group, with arithmetic subgroup $\Gamma$ and maximal compact subgroup $K$.  Let $X = \Gamma \backslash G/K$.  Let $\Delta$ be the Laplacian on $\Gamma\backslash G$, the image of the Casimir operator for the Lie algebra $\mathfrak{g}$; being $K$-invariant, the Laplacian descends to $X$. Let $\delta$ be the automorphic Dirac delta distribution centered at the basepoint $x_0 = \Gamma \cdot 1 \cdot K$ in $X$.  An automorphic heat kernel can be heuristically described as a function $u(x,t)$ on $X \times (0, \infty)$ satisfying the following two conditions:
$$(\partial_t \, - \, \Delta) \, u \; = \; 0 \hskip 1cm \text{and} \hskip 1cm \textrm{``} \; \displaystyle{\lim_{t \to 0^{+}}} \, u(x, t) \; = \; \delta(x) \; \textrm{''} \;.$$
We will give a precise definition of an automorphic heat kernel as a vector-valued function in \textsection\textsection\ref{more-precise}, after introducing the global automorphic Sobolev spaces where the automorphic heat kernel will take values.

\subsection*{Outline of the paper.} We begin, in  \textsection\ref{afc-sob-th}, by reviewing the results that we need from global automorphic Sobolev theory.  In  \textsection\ref{existence-sp-expn}, we give a precise definition of an automorphic heat kernel as a function taking values in global automorphic Sobolev spaces and construct it via its automorphic spectral expansion.  Note that this proves existence of the automorphic heat kernel without needing to discuss the convergence of the series obtained by winding up a heat kernel on a symmetric space.  The vector-valued (weak) integrals of Gelfand and Pettis are important here.  Then, in  \textsection\ref{uniqueness-section}, we prove the uniqueness of the automorphic heat kernel by proving the uniqueness of its automorphic spectral coefficient function.  Here we rely on the theory of strongly continuous semigroups of linear operators.  Finally, in  \textsection\ref{smoothness}, we prove the smoothness of the automorphic heat kernel as a function on $X$, for fixed $t > 0$, using a global automorphic Sobolev embedding theorem.  Appendices about differentiation and integration of vector-valued  functions, unbounded operators on Hilbert spaces, and strongly continuous semigroups are included for reference.

\section{Results from Global Automorphic Sobolev Theory}\label{afc-sob-th}

To state precisely and prove the automorphic spectral expansion for the automorphic heat kernel, we will use global automorphic Sobolev theory, developed in \cite{decelles2012}.  We briefly recall the results needed in this paper.  See \textsection2 (especially \textsection\textsection 2.3 and 2.4) of \cite{decelles2012} and Chapter 12 of \cite{garrett2018}.

First recall the spectral theory for $X = SL_2(\Z) \backslash \mathfrak{H}$.    For $f \in L^2(X)$, we have a spectral expansion in terms of eigenfunctions for the Laplacian $\Delta = y^2 (\tfrac{\partial^2}{\partial x^2} + \tfrac{\partial^2}{\partial y^2})$,
$$f \;\; \overset{L^2}{=} \;\; \sum_F \, \langle f, F\rangle \cdot F \; + \; \langle f, \Phi_0\rangle \cdot \Phi_0 \; + \; \frac{1}{4\pi i} \int_{\frac{1}{2} + i \R} \langle f, E_s \rangle \cdot E_s \, ds \;,$$
where $F$ ranges over an orthonormal basis of cusp forms, $\Phi_0$ is the constant automorphic form with unit $L^2$-norm, and $E_s$ is the real analytic Eisenstein series.  Note that the integrals (both the integral over $s$ and the integrals implied by the pairings) are not necessarily uniformly pointwise convergent, but they can be understood as extensions by continuity of continuous linear functionals on $C_c^\infty(X)$.

We can abbreviate (and generalize) this by denoting elements of the spectral ``basis'' (cusp forms, Eisenstein series, residues of Eisenstein series) uniformly as $\Phi_\xi$ for $\xi$ in a spectral parameter space $\Xi$ with Plancherel measure $d\xi$.  Then, for $f \in L^2(X)$, the spectral expansion is written
$$f \;\; = \;\; \int_\Xi \langle f, \Phi_\xi \rangle \cdot \Phi_\xi \; d\xi \; .$$
We can view $\Xi$ as a finite disjoint union of copies of spaces of the form $\Z^n \times \R^m$ for some $n, m \geq 0$, with the counting measure on each discrete space and the usual Euclidean measure on each Euclidean space.

Let $\lambda_\xi$ be the $\Delta$-eigenvalue of $\Phi_\xi$.  Since $\Delta$ is a nonpositive symmetric operator, $\lambda_\xi$ is nonpositive real.    Moreover, the map $\Lambda: \Xi \to \R$ by $\xi \mapsto \lambda_\xi$ is differentiable and of moderate growth, by a pre-trace formula.  (With suitable coordinates we expect it to be polynomial, up to powers of log, by Weyl's Law; see \cite{lapid-muller2009, matz2017}.)

For positive integer $s$, define an inner product $\langle \, , \, \rangle_s$ on $C_c^\infty(X)$ by
$$\langle \varphi,  \psi \rangle_s = \langle (1-\Delta)^s \varphi, \psi \rangle_{L^2} \; .$$
Let $H^s$ be the Hilbert space completion of $C_c^\infty(X)$ with respect to the topology induced by $\langle \, , \, \rangle_s$, and let $H^{-s}$ be its Hilbert space dual.  Note that $H^0 = L^2(X)$, and the Sobolev spaces form a nested family, $H^s \hookrightarrow H^{s - 1}$ for all $s$, with the inclusions being continuous linear maps.

The Laplacian acts nicely: $\Delta: H^s \to H^{s - 2}$ is a continuous linear map.  The spectral transform $\mathcal{F}: f \mapsto \langle f , \Phi_\xi \rangle$ gives a Hilbert space isomorphism to a weighted $L^2$-space $V^s$ on the spectral side:
$$V^s \;\; = \;\; \{f \text{ measurable} : (1-\lambda_\xi)^{\nicefrac{s}{2}} f(\xi) \in L^2(\Xi)\} \;. $$
These spaces form a nested family with continuous inclusions $V^s \hookrightarrow V^{s-1}$.

Denoting by $\mu$ the multiplication-by-$(1-\lambda_\xi)$ map, we have the following commutative diagram of Hilbert space isomorphisms:
\begin{equation*}
\xymatrix@C=2.2pc@R=3.2pc{
\dots  \;\;  H^{+s} \;\; \ar[r]^{(1-\Delta)}_\approx \ar[d]_{\mathcal{F}}^{\approx} & \;\; H^{+s-2}\;\;  \ar[r]^{(1-\Delta)}_\approx \ar[d]_{\mathcal{F}}^{\approx} & \;\;\dots \;\;  \ar[r]^{(1-\Delta)}_\approx & \;\; H^{-s+2}\;\;   \ar[r]^{(1-\Delta)}_\approx\ar[d]_{\mathcal{F}}^{\approx} & \;\; H^{-s}\;\;  \ar[d]_{\mathcal{F}}^{\approx}  \;\; \dots\\
\dots  \;\; V^{+s} \;\; \ar[r]^{\mu}_\approx & \;\; V^{+s-2} \;\; \ar[r]^{\mu}_\approx & \;\; \dots \;\; \ar[r]^{\mu}_\approx & \;\; V^{-s+2}  \ar[r]^{\mu}_\approx & \;\; V^{-s}  \;\; \dots }
\end{equation*}
This allows us to conclude that every $u \in H^s$ has a spectral expansion, converging in the $H^{s}$-topology:
$$u \;\; = \;\;\int_\Xi \mathcal{F}u(\xi) \cdot \Phi_\xi \; d\xi \;,$$
where the spectral coefficients $\mathcal{F}u$ lie in the weighted $L^2$-space $V^s$.

Moreover, a global automorphic Sobolev embedding theorem, which states that for $s > k + (\mathrm{dim}\, X)/2$, there is a continuous inclusion $H^s \hookrightarrow C^k$, allows comparison to $C^k$ convergence, and, in particular, uniform pointwise convergence, if desired.  Further, the global automorphic Sobolev embedding theorem implies that 
$$H^{\infty} \;\; = \;\; \bigcap_{s \in \Z} H^s \;\; \subset \;\; \bigcap_{k \in \Z} C^k \;\; = \;\; C^\infty \; ,$$
so an element lying in all global automorphic Sobolev spaces is in fact a smooth function on $X$.

Finally, we note that the automorphic Dirac delta distribution lies in every global automorphic Sobolev space of index \emph{strictly less} than $-(\mathrm{dim} \, X)/2$ and thus has the following automorphic spectral expansion
$$\delta \;\; = \;\; \int_\Xi \overline{\Phi}_\xi(x_0) \, \Phi_\xi \; d\xi \; ,$$
converging in the corresponding global automorphic Sobolev topologies.

\section{Existence and Spectral Expansion of an Automorphic Heat Kernel}\label{existence-sp-expn}

\subsection{An automorphic heat kernel as a $H^s$-valued function of $t$} \label{more-precise}

Having defined global automorphic Sobolev spaces, we can now make our notion of automorphic heat kernel more precise, using vector-valued functions.  Let $\ell$ be the smallest integer \emph{strictly greater} than $\mathrm{dim}\,  X/2$, so that $\delta$ lies in $H^{-\ell}$ but in no global automorphic Sobolev space of greater integer index.  We define an \emph{automorphic heat kernel} to be a map $U: (0, \infty) \to H^{-\ell}$ with the following properties.
\begin{enumerate}[label=(\roman*)]
\item $U$ satisfies the ``initial condition,'' $\displaystyle{\lim_{t \to 0^+} \, U(t) \; = \; \delta}$ in the topology of  $H^{-\ell}$.
\item For some $s \leq -\ell-2$, \, $U$ is strongly differentiable as a $H^s$-valued function and 
satisfies the ``heat equation,'' i.e. for $t >0$,
$$\tfrac{d}{dt} \,  U(t) \; -\;  \Delta  \, U(t) \;\; = \;\;0  \hskip .5cm \text{in } H^{s} .$$
\end{enumerate}

\begin{remark*} When we say that we differentiate $U$ as a $H^{s}$-valued function, we are really differentiating the function $\iota \circ  U: (0, \infty) \to H^{s}$, where $\iota$ is the inclusion $H^{-\ell} \hookrightarrow H^s$.  Moreover, when we say we view $\Delta U(t)$ as a $H^s$-valued function in the differential equation, we are really considering the function $j \circ \Delta U(t)$, where $j$ is the inclusion $H^{-\ell-2} \hookrightarrow H^{s}$. Both of these inclusions, $\iota$ and $j$, are continuous and linear, so for brevity, we suppress them from the notation.
\end{remark*}

\subsection{Translation to the spectral side}

Our strategy is to apply a spectral transform, solve a simpler differential equation on the spectral side, and then apply spectral inversion to obtain our solution on the physical side.  The equivalences in the following lemma assert the validity of such translation between the physical and spectral sides.  In particular, part \ref{trans-heat-eqtn} tells us that the $H^s$-valued function $U$ is a solution to the automorphic heat equation if and only if $\mathcal{F} \circ U$ is a solution to the differential equation of $V^{s}$-valued functions: $Y'(t) = (1-\mu)Y(t)$.

\begin{lemma}[Translation Lemma] \label{translation} Consider a vector-valued function $U:(0, \infty) \to H^s$ and an element $\theta \in H^s$.  Let $\iota_H$ and $\iota_V$ denote the continuous inclusions $H^s \hookrightarrow H^{s-2}$ and $V^s \hookrightarrow V^{s-2}$ respectively.
\begin{enumerate}[label=(\roman*)]
\item\label{limit} $\displaystyle{\lim_{t\to 0^+} U(t) \; = \; \theta}\;\;$ in the $H^{s}$-topology  if and only if \, $\displaystyle{\lim_{t\to 0^+} (\mathcal{F}\circ U)(t) \; = \; \mathcal{F}\theta}\;\;$ in the $V^{s}$-topology.
\item\label{differentiability} $U$ is weakly (or strongly) differentiable (or $C^1$) if and only if $\mathcal{F}\circ U$ is weakly (resp. strongly) differentiable (resp. $C^1$).
\item\label{trans-heat-eqtn} As $H^{s-2}$-valued maps, $\iota_H \circ U' \; = \; \Delta \circ U$  if and only if\, $\iota_V \circ (\mathcal{F} \circ U)' \; = \; (1-\mu)\circ (\mathcal{F} \circ U)$, as $V^{s-2}$-valued maps.
\end{enumerate}
Here, we use the prime notation to refer both to weak and strong derivatives. 
\begin{proof}  Parts \ref{limit} and \ref{differentiability} follow immediately from the continuity and linearity of $\mathcal{F}$ and $\mathcal{F}^{-1}$ and Lemma \ref{comp-w-diffble} in Appendix \ref{lim-deriv-V-vald}.  To prove part \ref{trans-heat-eqtn}, we apply $\mathcal{F}$ to both sides of $\iota_H \circ U'  \, = \, \Delta \circ U$.  Note that $\mathcal{F} \circ \iota_H = \iota_V \circ \mathcal{F}$.  Then since $\mathcal{F}$ is continuous and linear, we may use Lemma \ref{comp-w-diffble} again, and we have
$$\mathcal{F}\circ \iota_H \circ U' \;\; = \;\; \iota_V \circ \mathcal{F} \circ U' \;\; = \;\; \iota_V \circ (\mathcal{F} \circ U)' \;.$$
On the other hand, from \cite{decelles2012}, we know $\mathcal{F} \circ \Delta = (1-\mu)\circ \mathcal{F}$.  Thus, applying $\mathcal{F}$ to both sides of $\iota_H \circ U' \; = \; \Delta \circ U$ yields $\iota_V \circ (\mathcal{F} \circ U)' \; = \; (1-\mu)\circ (\mathcal{F} \circ U)$.  Since $\mathcal{F}$ is an isomorphism this means that the two equations are equivalent. \end{proof}
\end{lemma}

\subsection{Construction of automorphic heat kernel via spectral synthesis}
\label{spectral-constr}

Next we define a (function-on-$\Xi$)-valued map, which, for fixed ``time'' $t$, gives the expected automorphic spectral coefficient, as a function of $\xi$, for the automorphic heat kernel: for $t \geq 0$, we define $\widetilde{U}(t)$ to be the map $\Xi \to \C$ given by 
$$\widetilde{U}(t): \; \xi \; \mapsto \; \overline{\Phi}_\xi(x_0) \, e^{\lambda_\xi \, t}.$$
Since $\lambda_\xi$ is real and nonpositive, we often write $e^{\lambda_\xi \, t}$ as $e^{-|\lambda_\xi| \, t}$, to aid our intuition. 

\begin{prop}\label{spectral-V-valued-initial-cond} The (function-on-$\Xi$)-valued map $\widetilde{U}$ defined above has the following two properties.
 \begin{enumerate}[label=(\roman*)] \item \label{V-valued} For $t \geq 0$, $\widetilde{U}(t) \in V^{-\ell}$.

\item \label{spectral-initial-cond}  $\displaystyle{\lim_{t\to 0^+} \; \widetilde{U}(t) \; = \; \mathcal{F}\delta\;}$ in the topology of $V^{-\ell}$.
\end{enumerate}
\begin{proof}  We define $\mathcal{E}: (0, \infty) \to L^\infty(\Xi)$ by $\mathcal{E}(t): \xi \, \mapsto \, e^{-|\lambda_\xi| \, t}$.  Since $\ell$ was chosen to ensure $\delta \in H^{-\ell}$, we have $\mathcal{F}\delta \in V^{-\ell}$.  Then, by the H\"{o}lder inequality, $\widetilde{U}(t) = \mathcal{F}\delta \cdot \mathcal{E}(t)$ is also in $V^{-\ell}$, proving \ref{V-valued}.

To show that $\widetilde{U}(t) \to \mathcal{F}\delta$ as $t \to 0^+$ in the $V^{-\ell}$ topology, we show that $\widetilde{U}(t)$ approaches $\mathcal{F}\delta$ weakly and in the norm.  Recall that the continuous linear dual of $V^{-s}$ is $\big(V^{-s}\big)^\ast = V^{s}$.  Let $\langle \, , \, \rangle$ denote the pairing on $V^{s} \times V^{-s}$.  Take any $f \in V^{\ell}$. Lebesgue Dominated Convergence implies
$$\lim_{t\to 0^+} \; \int_\Xi f(\xi) \overline{\Phi}_\xi(x_0) e^{-|\lambda_\xi| \, t} \, d\xi  \;\; = \;\; \int_\Xi f(\xi) \overline{\Phi}_\xi(x_0) \, d\xi \;,$$
i.e. $\langle f, \widetilde{U}(t) \rangle \to \langle f, \mathcal{F}\delta\rangle$ as $t \to 0^+$.  Moreover, by Lebesgue Dominated Convergence,
$$\lim_{t \to 0^+} \; \int_\Xi (1+|\lambda_\xi|)^{-\ell} e^{-2|\lambda_\xi| t} \, |\overline{\Phi}_\xi(x_0)|^2 \, d\xi \;\; = \;\;  \int_\Xi (1+|\lambda_\xi|)^{-\ell}  \, |\overline{\Phi}_\xi(x_0)|^2  \, d\xi\;, $$
i.e. $\lVert \widetilde{U}(t)\rVert_{V^{-\ell}}^2$ approaches $\lVert \mathcal{F}\delta\rVert_{V^{-\ell}}^2$ as $t\to 0^+$, completing the proof of \ref{spectral-initial-cond}. \end{proof}
\end{prop}

\begin{prop}\label{spectral-diff-eqtn}  When viewed as a $V^{-\ell-3}$-valued map, $\widetilde{U}$ is weakly continuously differentiable on $(0, \infty)$, and its weak derivative $\left(\tfrac{d}{dt}\right)^w  \,\widetilde{U}$ satisfies:
$$\left(\tfrac{d}{dt}\right)^w  \,\widetilde{U}(t) \; - \; (1-\mu) \widetilde{U}(t) \;\; = \;\; 0 \hskip 1cm \text{ in } V^{-\ell-3} \;,$$
where, as above, $(1-\mu)$ is the multiplication-by-$\lambda_\xi$ map.  Moreover, when viewed as a $V^{-\ell-5}$-valued map, $\widetilde{U}$ is strongly continuously differentiable on $(0, \infty)$, and its strong derivative $\tfrac{d}{dt}  \,\widetilde{U}$ satisfies
$$\tfrac{d}{dt}  \,\widetilde{U}(t) \; - \; (1-\mu) \widetilde{U}(t) \;\; = \;\; 0 \hskip 1cm \text{ in } V^{-\ell-5} \;.$$

\begin{remark*} When we say that we view $\widetilde{U}$ as a $V^{-\ell-3}$-valued (or $V^{-\ell-5}$-valued) map, we are really considering the map $\iota \circ \widetilde{U}: (0, \infty) \to V^{-\ell-3}$, where $\iota$ is the inclusion $V^{-\ell} \hookrightarrow V^{-\ell-3}$ (resp. $V^{-\ell} \hookrightarrow V^{-\ell-5}$).  Since each of these inclusions is continuous and linear, Lemma \ref{comp-w-diffble} in Appendix \ref{lim-deriv-V-vald} ensures that we do not lose anything by suppressing them from the notation.
\end{remark*}

\begin{proof} We show that $\widetilde{U}$ is weakly $C^k$, when viewed as a $V^{-\ell-2N}$-valued function, for $N > k$.  In particular, this will show that it is weakly $C^1$, when viewed as a $V^{-\ell-3}$-valued function.  To draw conclusions about strong differentiability, we apply the weak-to-strong differentiability principle: weak $C^k$ implies strong $C^{k-1}$; see Appendix \ref{lim-deriv-V-vald}.  This implies that $\widetilde{U}$ is strongly $C^{k-1}$ when viewed as a $V^{-\ell -2N}$-valued function for $N>k>1$, so in particular, it is strongly $C^1$ when viewed as a $V^{-\ell-5}$-valued function.

Fix $k \geq 1$, and let $s = \ell+2N$ for some $N>k$.  Recall that the continuous linear dual of $V^{-s}$ is $\big(V^{-s}\big)^\ast = V^{s}$.  Let $\langle \, , \, \rangle$ denote the pairing on $V^{s} \times V^{-s}$.  
Take any $f \in V^{s}$.  We will show that 
$$t \;\;\; \mapsto \;\;\; \langle f, \widetilde{U}(t)\rangle \;\; = \;\; \int_\Xi f(\xi) \overline{\Phi}_\xi(x_0) e^{-|\lambda_\xi| t} \, d\xi \;.$$
is a $C^k$ scalar-valued function.  The key is to identify this function with a \emph{Gelfand-Pettis integral} of a $C^k(0,\infty)$-valued function of $\xi$.  (See Appendix \ref{gelfand-pettis}.)

Consider the continuous vector-valued function $\mathcal{W}: \Xi \to C^k(0, \infty)$ given by $\mathcal{W}(\xi): t \mapsto (1-\lambda_\xi)^{-N}\, e^{-|\lambda_\xi| t}$.   Since $\mathcal{W}$ is not compactly supported on $\Xi$, we need to compactify $\Xi$ and give it a finite measure, in order to ensure that the Gelfand-Pettis integral of $\mathcal{W}$ exists.

We view $\Xi$ is a finite disjoint union of copies of spaces of the form $\Z^n \times \R^m$ for some $n, m \geq 0$, with the usual metrics on $\Z^n$ and $\R^m$.  Let $\widehat{\Xi}$ be the compactification of $\Xi$ obtained by taking the disjoint union of all the Cartesian products $\widehat{(\Z^n)} \times \widehat{(\R^m)}$, where the hat denotes the usual one-point compactification of $\Z^n$ or $\R^m$.   We adjust $d\xi$ to obtain a finite measure for $\widehat{\Xi}$, by 
$$d\hat{\xi} \;\; = \;\; (1-\lambda_\xi)^N \, f(\xi) \, \overline{\Phi}_\xi(x_0) \, d\xi \;.$$
Since $f \in V^{\ell + 2N}$ and $\overline{\Phi}_\xi(x_0) \in V^{-\ell}$, the ``adjustment'' function, $(1-\lambda_\xi)^N  f(\xi) \,\overline{\Phi}_\xi(x_0)$, is in $L^1(\Xi)$, i.e. $d\hat{\xi}$ is a finite measure on $\Xi$ and thus also on $\widehat{\Xi}$, since $\widehat{\Xi}$ differs from $\Xi$ only by a finite set of points.

We extend $\mathcal{W}: \Xi \to C^k(0, \infty)$ to $\widehat{\Xi}$, by the zero function, which is certainly in $C^k(0, \infty)$.  The extension $\mathcal{W}: \widehat{\Xi} \to C^k(0, \infty)$ is continuous, since, for any compact $C \subset (0, \infty)$ and any nonnegative $n \leq k$, 
$$\sup_{t \in C} \big| \partial_t^n (1-\lambda_\xi)^{-N} e^{-|\lambda_\xi| t} \big| \, = \,  |\lambda_\xi|^n (1+|\lambda_\xi|)^{-N} \sup_{t \in C} e^{-|\lambda_\xi| t}\, \leq  \, |\lambda_\xi|^k \, (1+|\lambda_\xi|)^{-N}, $$
which approaches zero as $|\lambda_\xi| \to \infty$, since $N > k$.

The extension of $\mathcal{W}$ to $\widehat{\Xi}$, which we also denote $\mathcal{W}$, is continuous on $\widehat{\Xi}$, which is a \emph{compact} Hausdorff topological space with a \emph{finite} positive regular Borel measure $d\hat{\xi}$, and $\mathcal{W}$ takes values in the quasi-complete, locally convex topological vector space $C^k(0, \infty)$.  Therefore $\mathcal{W}$ has a $C^k(0, \infty)$-valued Gelfand-Pettis integral.  Since evaluation-at-$t$ is a continuous linear functional, call it $E_t$, on $C^k(0, \infty)$, it commutes with the Gelfand-Pettis integral of $\mathcal{W}$, and we have:
$$ \hskip -2cm E_t\bigg(\int_{\widehat{\Xi}} \mathcal{W}(\xi) \, d \hat{\xi}\bigg) \;\;= \;\; \int_{\widehat{\Xi}} E_t\big(\mathcal{W}(\xi)\big) \, d \hat{\xi} \;\; = \;\; \int_{\widehat{\Xi}}(1-\lambda_\xi)^{-N} \, e^{-|\lambda_\xi| t}\, d \hat{\xi}$$

$$ \;\; = \;\; \int_{\widehat{\Xi}}(1-\lambda_\xi)^{-N} \, e^{-|\lambda_\xi| t}\, (1-\lambda_\xi)^N f(\xi) \, \overline{\Phi}_\xi(x_0) \; d\xi \;\; = \;\; \int_{\widehat{\Xi}} f(\xi) \, \overline{\Phi}_\xi(x_0) \, e^{-|\lambda_\xi| t} \; d\xi \;.$$
Since $\Xi$ and $\widehat{\Xi}$ differ only by a finite set, 
$$\int_{\widehat{\Xi}} f(\xi) \, \overline{\Phi}_\xi(x_0) \, e^{-|\lambda_\xi| t} \; d\xi  \;\; = \;\; \int_\Xi f(\xi) \overline{\Phi}_\xi(x_0) e^{-|\lambda_\xi| t} \, d\xi \;\; = \;\; \langle f, \widetilde{U}(t) \rangle \;.$$
Thus the map $t \mapsto \langle f, \widetilde{U}(t) \rangle$ is precisely the $C^k(0, \infty)$-valued Gelfand-Pettis integral of $\mathcal{W}$ over $\widehat{\Xi}$.  Since this holds for any $f \in V^{\ell+2N}$, this completes the proof that $\widetilde{U}$ is weakly $C^k$, when considered as a $V^{-\ell-2N}$-valued function, for $N>k$.

Moreover, since $\tfrac{d}{dt}$ is a continuous linear operator on $C^k(0, \infty)$, it also commutes with the Gelfand-Pettis integral of $\mathcal{W}$.  So for any $f \in V^{\ell+2N}$,
$$\tfrac{d}{dt}(t \mapsto \langle f, \widetilde{U}(t) \rangle) \;\; = \;\; \tfrac{d}{dt} \bigg(\int_{\widehat{\Xi}} \mathcal{W}(\xi) \, d \hat{\xi}\bigg) \;\;= \;\; \int_{\widehat{\Xi}} \tfrac{d}{dt}\big(\mathcal{W}(\xi)\big) \, d \hat{\xi} \;,$$
and the derivative is easily computed as
$$\tfrac{d}{dt}\big( \mathcal{W}(\xi) \big) \;\; = \;\; \tfrac{d}{dt} \bigg( t \mapsto (1-\lambda_\xi)^{-N} \, e^{\lambda_\xi t}\bigg) \;\; = \;\; t \mapsto \lambda_\xi (1-\lambda_\xi)^{-N} \, e^{\lambda_\xi t} \;\; = \;\; \lambda_\xi \cdot \mathcal{W}(\xi) \; ,$$
so 
$$\int_{\widehat{\Xi}} \tfrac{d}{dt}\big(\mathcal{W}(\xi)\big) \, d \hat{\xi} \;\; = \;\; \int_{\widehat{\Xi}} \lambda_\xi \cdot \mathcal{W}(\xi) \, d\hat{\xi} \;\; = \;\; \bigg( t \mapsto \int_\Xi f(\xi) \cdot \lambda_\xi \, \overline{\Phi}_\xi(x_0) \, e^{\lambda_\xi t} \; d\xi \bigg) \;.$$
On the other hand,
$$\langle f, (1-\mu)\widetilde{U}(t) \rangle \;\; = \;\; \int_\Xi f(\xi) \; \lambda_\xi \cdot  \overline{\Phi}_\xi(x_0) \, e^{\lambda_\xi t} \; d\xi \;.$$
Thus, for all $f \in V^{\ell+2N}$, 
$$\tfrac{d}{dt} \big(t \mapsto \langle f, \widetilde{U}(t) \rangle\big) \;\; = \;\; t\mapsto \langle f, (1-\mu) \, \widetilde{U}(t)\rangle \;, $$
i.e. $(1-\mu)\widetilde{U}(t)$ is the weak derivative of $\widetilde{U}(t)$, when $\widetilde{U}$ is considered as a $V^{-\ell-2N}$-valued function.

Now suppose $k>1$.  In this case, as discussed above, $\widetilde{U}$ is also strongly differentiable, in fact strongly $C^{k-1}$.  Moreover the strong derivative is equal to the weak derivative; thus we have $$\tfrac{d}{dt} \widetilde{U}(t) \; - \; (1-\mu) \widetilde{U}(t) \;\; = \;\; 0 \;,$$
which completes the proof of the proposition.\end{proof}
\end{prop}

\begin{remark}\label{spec-coef-props} We have shown that for $s < -\ell-4$, $\widetilde{U}$ is strongly continuous as a $V^s$-valued function on $[0, \infty)$, is strongly differentiable as a $V^s$-valued function on $(0, \infty)$, takes values in $V^{s+2}$, and solves the initial value problem:
$$\tfrac{d}{dt} Y(t) \;\; = \;\; A \; Y(t) \; ; \hskip .5cm Y(0) \; = \;\; \mathcal{F}\delta\;.$$
Later, we will use semigroup theory to show that $\widetilde{U}$ is the only such solution to this initial value problem.  See Propositions \ref{uniqueness-coef} and \ref{uniqueness-abstract}.
\end{remark}

We are now ready to prove our first major result, the construction of an automorphic heat kernel via spectral synthesis.

\begin{theorem}\label{constr-afc-heat-ker} For $t \geq 0$, let $U(t) = (\mathcal{F}^{-1} \circ \widetilde{U})(t)$, i.e. we define $U(t)$ by its automorphic spectral expansion,   
$$U(t) \;\; = \;\; \int_\Xi \overline{\Phi}_\xi(x_0) \cdot e^{\lambda_\xi\,  t} \cdot \Phi_\xi \, d\xi \;,$$
converging in the $H^{-\ell}$-topology, where, as in the above discussion, $\ell$ is the smallest integer strictly greater than $\mathrm{dim}\,  X/2$, so that $\delta$ lies in $H^{-\ell}$ but not $H^{-\ell+1}$.  

 \begin{enumerate}[label=(\roman*)] \item \label{H-valued} For $t \geq 0$, $U(t) \in H^{-\ell}$.
\item \label{initial-cond} $\displaystyle{\lim_{t\to0^+} U(t) \;=\; \delta}$\, in the topology of $H^{-\ell}$.
\item For $s < -\ell-2$, viewing $U$ as a $H^s$-valued function, $U$ is weakly continuously differentiable on $(0, \infty)$ and satisfies the ``weak heat equation,'' i.e. for $t >0$,
$$\left(\tfrac{d}{dt}\right)^w U(t) \; - \; \Delta \, U(t) \;\; = \;\; 0 \hskip .5cm \text{ in } H^{s} ,$$
where $\left(\tfrac{d}{dt}\right)^w U$ denotes the weak derivative of $U$.
\item \label{ker-strong-diffble}  For $s < -\ell -4$, viewing $U$ as a $H^s$-valued function, $U$ is strongly continuously differentiable on $(0, \infty)$ and 
satisfies the ``heat equation,'' i.e. for $t >0$,
$$\tfrac{d}{dt} U(t) \; -\;  \Delta  \, U(t) \;\; = \;\;0  \hskip .5cm \text{in } H^{s},$$
where $\tfrac{d}{dt} U$ denotes the strong derivative of $U$.
\end{enumerate}
Thus $U(t)$ is an automorphic heat kernel, according to the definition in \textsection\textsection\ref{more-precise}.
\begin{proof}  These results follow immediately from the translation lemma (Lemma \ref{translation}) and the corresponding results (Propositions \ref{spectral-V-valued-initial-cond} and \ref{spectral-diff-eqtn}) on the spectral side. \end{proof}
\end{theorem}


\section{Uniqueness of the Automorphic Heat Kernel} \label{uniqueness-section}

We will prove the uniqueness of the automorphic heat kernel $U$ by proving the uniqueness of its automorphic spectral coefficient function, $\widetilde{U}$.  In particular, we will prove the uniqueness of (suitable) solutions to the initial value problem
$$\tfrac{d}{dt} Y(t) \;\; = \;\; (1-\mu) \; Y(t) \; ; \hskip .5cm  Y(0) = \mathcal{F} \delta \;,$$
where, as above, the map $(1-\mu)$ is multiplication by $\lambda_\xi$ and $\mathcal{F} \delta$ is the automorphic spectral coefficient function for the automorphic delta function at the basepoint: $\mathcal{F}\delta(\xi) = \overline{\Phi}_\xi(x_0)$. Note that $(1-\mu)$ is a continuous linear map when considered as a map from $V^{-\ell}$  to $V^{-\ell-2}$ but \emph{not} when considered as a map from $\iota(V^{-\ell})$ to $V^{-\ell-2}$, where $\iota$ is the continuous inclusion $V^{-\ell} \hookrightarrow V^{-\ell-2}$.  Thus we \emph{cannot} use results about differential equations of vector-valued functions of the form  $\tfrac{d}{dt} Y(t) \, = \, B Y(t)$, where $B$ is a bounded linear operator on a vector space.  Instead, we consider $(1-\mu)$ as an unbounded operator, $M$, on $V^{-\ell-2}$ with domain $V^{-\ell}$ and use semigroup theory to show that the initial value problem has a unique solution.  See Appendices \ref{unbdd-ops} and \ref{semigroup-theory} for background on unbounded operators on Hilbert spaces and semigroup theory, respectively.

\subsection{Multiplication Operators on Weighted $L^2$-spaces} \label{mult-ops}

This discussion generalizes and adapts Grubb's discussion in \cite{grubb2009} of multiplication operators on $L^2(\Omega)$, where $\Omega$ is an open subset of $\R^n$.  Let $V=V^s$ for some $s$, and let $M$ be the unbounded operator on $V$ given by 
$$\big(Mv\big)(\xi) \;\; = \;\; \lambda_\xi \cdot v(\xi)$$
with domain $D \; = \; \mathrm{Dom} \, M \; = \; \{ v \in V : \, Mv \in V \}$.

\begin{prop}\label{domain} The domain $D$ of $M$ is equal to $V^{s+2}$.
\begin{proof}  By definition, $V^{s+2} =  \{v \textrm{ measurable }: (I-M)v \, \in \, V\}$. For $v \in D$, we have $v$ and $Mv$ in $V$ and thus $(I-M)v$ is as well; thus $D \subset V^{s+2}$.  Now take any $v \in V^{s+2}$, i.e. $v$ is measurable and $(I-M) v \, \in \, V$.  To show $v \in D$, it suffices to show $v \in V$, since this would imply $Mv \in V$, by closure.  Since $\lambda_\xi$ is real and nonpositive, 
$\big| \, (1-\lambda_\xi) \, v(\xi) \, \big| \; \geq \; |v(\xi)| \;.$  Thus the fact that $(I-M)v$ is in $V^s=V$ implies that $v$ is also.  Thus $V^{s+2} \subset D$. \end{proof}
\end{prop}

\begin{prop} \label{dd-self-adjt} $M$ is a densely defined self-adjoint operator on $V$. 
\begin{proof}  This proof is an adaptation of the proof of Theorem 12.13 in \cite{grubb2009}.  Let $\langle \, , \, \rangle$ be the inner product on $V = V^{s}$, i.e. for $v, w \in V$,
$$\langle v, w \rangle \;\; = \;\; \int_{\Xi} v(\xi) \cdot \overline{w}(\xi) \, (1-\lambda_\xi)^{s} \, d\xi \; .$$
Let $\Lambda: \Xi \to \C$ denote the function $\xi \mapsto \lambda_\xi$.  We observe that a function lies in the domain of $M$ if and only if it is a quotient of the form $v/(1-\Lambda)$ for some $v \in V$.  Note also that if $w \in V$, then $w/(1-\Lambda) \in V$ since 
$$ \left| \frac{w(\xi)}{1-\lambda_\xi} \right| \;\; = \;\; \frac{|w(\xi)|}{1+|\lambda_\xi|} \;\; \leq \;\;|w(\xi)| \;.$$
The domain of $M$ is certainly dense in $V$ for if not, then there would be a nonzero $w \in V$ with $w \perp \overline{D}$, i.e.
$$0 \;\; = \;\; \left\langle \frac{v}{1-\Lambda} , \, w \right\rangle \;\; = \;\; \left\langle v, \, \frac{w}{1-\Lambda}\right\rangle \hskip 1cm \textrm{for all } v \in V \; ,$$
implying $w/(1-\Lambda) =0$, contradicting $w \neq 0$.
 
Clearly $M$ is symmetric since $\lambda_\xi$ is real.  Since $M$ is densely defined, it has a unique adjoint $M^\ast$, whose domain is
$$\mathrm{Dom} \, M^\ast  = \{w \in V : \textrm{there is } w_1 \in V \textrm{ with } \langle Mv, w \rangle = \langle v, w_1\rangle \textrm{ for all } v \in D \} \,,$$
by Remark \ref{dom-adjt}.  Since $M$ is densely defined and symmetric, to show that $M$ is self-adjoint, it suffices to show $\mathrm{Dom} \, M^\ast  \subset  D$, by Remark \ref{dd-sym-self-adjt}.  Take $w \in \mathrm{Dom} \, M^\ast$, and let $w_1$ satisfy
$$\langle Mv, w \rangle \;\; =\;\; \langle v, w_1 \rangle \hskip 1cm \textrm{ for all } \;v \in D \;.$$
Since any element of $D$ is a quotient of the form $u/(1-\Lambda)$, where $u \in V$, we can rewrite this as
$$\left\langle \frac{\Lambda \, u}{1-\Lambda}, w \right\rangle \;\; = \;\; \left\langle \frac{u}{1-\Lambda}, \, w_1 \right\rangle \hskip 1cm \textrm{for all } \; u \in V \; .$$
As noted above, $w_1/(1-\Lambda) \in V$ since $w_1$ is.  By a similar argument $\Lambda w/(1-\Lambda)$ is also in $V$ since $w$ is.  Therefore we have
$$\left\langle u , \, \frac{\Lambda \, w}{1-\Lambda} \right\rangle \;\; = \;\; \left\langle u, \, \frac{w_1}{1-\Lambda} \right\rangle \hskip 1cm \textrm{for all } \; u \in V \; ,$$
which implies $Mw = w_1$.  Since $w_1 \in V$, we have $Mw \in V$ and thus $w \in D$.  We have shown $\mathrm{Dom} \, M^\ast \, \subset \, D$, completing the proof that $M$ is self-adjoint.
\end{proof}\end{prop}

\begin{prop}\label{neg-res} $M$ is negative, and its resolvent set contains $(0, \infty)$.
\begin{proof} The negativity of $M$ follows easily from the nonpositivity of $\lambda_\xi$.  Recall that the resolvent set of $M$ consists of all $C \in \C$ such that the operator $M-CI$ is a bijection $D \to V$ and $(M-CI)^{-1}$ is bounded.  Take any positive real number $C$.  To show injectivity, consider $v \in D$ in the kernel of $M-CI$.  Then 
$$0 \;\; =\;\; (M-CI)v \,(\xi) \;\; = \;\; -(|\lambda_\xi|+C)\,  v(\xi) \hskip 1cm \textrm{a.a. } \xi \in \Xi \;,$$
implying $v = 0$ in $V$.  To show surjectivity, take $w \in V$, and let $v = w/(\Lambda-C)$.  Certainly $(M-CI)v = w$, as long as $v$ is in the domain of $M$.  As shown in the proof of Proposition \ref{dd-self-adjt}, to see that $v \in D$, it suffices to show $v = u/(1-\Lambda)$ for some $u \in V$.  Taking $u = ((1-\Lambda)/(\Lambda-C))w$ suffices.  Finally, the resolvent $(M-CI)^{-1}$ is bounded, since, for any $w \in V$, $\lVert (M-CI)^{-1}w \rVert   \leq  C^{-1} \lVert w \rVert$.
\end{proof}
\end{prop}

\subsection{Uniqueness and an improvement on the continuous differentiability of the automorphic heat kernel} \label{uniqueness-afc-heat-kernel}

We have shown that the map $(1-\mu)$, considered as an unbounded operator on $V^{s}$, is densely-defined, self-adjoint, negative, and has the positive real numbers in its resolvent set.  Semigroup theory implies the existence and uniqueness of a suitable solution to the initial value problem:
$$\tfrac{d}{dt} Y(t) \;\; = \;\; (1-\mu) \; Y(t) \; ; \hskip .5cm  Y(0) = \mathcal{F} \delta \;,$$
as we will explicate and prove below.  (See also  \textsection\ref{acp-uniqueness}, on the uniqueness of solutions of initial value problems of this form, in the appendix on semigroup theory.)

Besides allowing us to prove the uniqueness of the spectral coefficient function $\widetilde{U}$ for the automorphic heat kernel, this will also allow us to improve Proposition \ref{spectral-diff-eqtn}, by showing that $\widetilde{U}$ is strongly continuously differentiable as a $V^{-\ell-2}$-valued function on $[0, \infty)$, not merely as a $V^{-\ell-5}$-valued function on $(0, \infty)$.   See Proposition \ref{uniqueness-coef}.

Moreover, translating these results to the physical side will give our second major result, Theorem \ref{uniqeness-afc-heat-ker}: the automorphic heat kernel $U(t)$ constructed in Theorem \ref{constr-afc-heat-ker} via its automorphic spectral expansion is the \emph{unique} automorphic heat kernel, as defined precisely in  \textsection\ref{more-precise}, and, when viewed as a $H^{-\ell-2}$-valued function, is strongly continuously differentiable on $[0, \infty)$. 

\begin{prop}\label{uniqueness-coef} $\widetilde{U}(t)$ is the unique solution to the initial value problem, 
\begin{equation}\label{IVP} 
\tfrac{d}{dt} Y(t) \; = \; (1-\mu) Y(t)\; ; \hskip 1cm Y(0) \; = \; \mathcal{F} \delta \;,
\tag{$\ast$} \end{equation}
in the following sense: for any $s \leq -\ell-2$, there is no other function besides $\widetilde{U}$ that (i) is strongly continuous as a $V^s$-valued function on $[0, \infty)$, (ii) is strongly differentiable as a $V^s$-valued function on $(0, \infty)$, (iii) takes values in $V^{s+2}$, and (iv) solves the initial value problem (\ref{IVP}) where the differentiation is understood as strong differentiation of $V^s$-valued functions.  Moreover, $\widetilde{U}$ is in fact strongly continuously differentiable as a $V^{-\ell-2}$-valued function on $[0, \infty)$.
\begin{proof} Since, in this discussion, all continuity, differentiability, and continuous differentiability of vector-valued functions is strong, we simply write ``continuous'' for strongly continuous, etc.

For any fixed $s\in \Z$, we can view $(1-\mu)$ as an unbounded multiplication operator, $M$, on $V^{s}$ with domain $V^{s+2}$, by Proposition \ref{domain}.  By Propositions \ref{dd-self-adjt} and \ref{neg-res}, $M$ is densely defined, self-adjoint, negative, and has a resolvent set containing the positive real numbers.  Semigroup theory (see e.g. Corollary 14.11 and Theorem 14.2 in \cite{grubb2009}, restated below in Appendix \ref{semigroup-theory}, and Remark \ref{Grubb-14.2-C1}) then implies that $M$ is the infinitesimal generator for a strongly continuous contraction semigroup $G(t)$, and for any $v_0 \in V^{s+2}$,  the vector-valued function $t \mapsto G(t) v_0$ is continuously $V^s$-differentiable, takes values in $V^{s+2}$, and solves the initial value problem:
$$\tfrac{d}{dt} Y(t) \; = \; M \, Y(t)\; ; \hskip 1cm Y(0) \; = \; v_0.$$
Proposition \ref{uniqueness-abstract} implies that there is no other solution that is continuous as a $V^s$-valued function on $[0,\infty)$, differentiable as a $V^s$-valued function on $(0, \infty)$, and takes values in $V^{s+2}$.  In particular, since $\mathcal{F} \delta \, \in \, V^{-\ell}$, this construction yields a $V^{-\ell}$-valued solution, namely $t \mapsto G(t) \mathcal{F}\delta$, to the initial value problem (\ref{IVP}) that is continuously differentiable as a $V^{-\ell-2}$-valued function.

On the other hand, by Propositions \ref{spectral-V-valued-initial-cond} and \ref{spectral-diff-eqtn}, $\widetilde{U}$ is continuous on $[0, \infty)$ and continuously differentiable on $(0, \infty)$ when viewed as a $V^{-\ell-5}$-valued function, takes values in $V^{-\ell} \subset V^{-\ell-3}$, and  solves the ``same'' initial value problem, although differentiation must be understood as differentiation of a $V^{-\ell-5}$-valued function. (See Remark \ref{spec-coef-props}.)  We claim that $\widetilde{U}(t)$ is in fact the same as the solution $t \mapsto G(t) \mathcal{F} \delta$ coming from semigroup theory, and thus is in fact continuously differentiable when viewed as a $V^{-\ell-2}$-valued function on $[0, \infty)$.  The key point is that both functions, $t \mapsto G(t) \mathcal{F} \delta$ and $\widetilde{U}$, can be viewed as continuous $V^{-\ell-5}$-valued functions on $[0, \infty)$, differentiable $V^{-\ell-5}$-valued functions on $(0, \infty)$, that take values in $V^{-\ell-3}$ and solve the same initial value problem, so they must be equal by the uniqueness of such solutions (Proposition \ref{uniqueness-abstract}).  We prove this carefully in the following, slightly more general discussion.

For $s' \leq s$, let $M_s$ denote the multiplication operator from $V^{s+2} \to V^s$, and $M_{s'}$ the multiplication operator from $V^{s'+2} \to V^{s'}$.  Fix $v_0 \in V^{s+2}$, and let $Y_s(t)$ be the unique function that is continuous as a $V^s$-valued function on $[0, \infty)$, differentiable as a $V^s$-valued function on $(0, \infty)$, takes values in $V^{s+2}$, and solves 
\begin{equation}\label{s-IVP}
\tfrac{d}{dt} \; \iota \circ Y \; = \; M_s \circ Y\; ; \hskip 1cm Y(0) \; = \; v_0 \,, \tag{IVP--$s$}
\end{equation}
where $\iota$ is the continuous inclusion $V^{s+2} \hookrightarrow V^s$.  This emphasizes that the implied limit on the left side is taken with respect to the topology on $V^s$ rather than $V^{s+2}$.

Since $v_0 \in V^{s+2} \subset V^{s'+2}$, we may define $Y_{s'}$ as the unique function that is continuous as a $V^{s'}$-valued function on $[0, \infty)$, differentiable as a $V^{s'}$-valued function on $(0, \infty)$, takes values in $V^{s'+2}$, and solves the initial value problem
\begin{equation}\label{sprime-IVP}
\tfrac{d}{dt} \; \iota'\circ Y \; = \; M_{s'} \circ Y \; ; \hskip 1cm Y(0) \; = \; v_0 \, , \tag{IVP--$s'$}\end{equation}
where $\iota'$ is the continuous inclusion $V^{s'+2} \hookrightarrow V^{s'}$.  We check that $Y_{s'}$ is the same as $Y_s$.  This follows from the fact that $Y_s$ satisfies (\ref{sprime-IVP}), which we now show.

Let $j_0$ be the continuous inclusion $V^s \hookrightarrow V^{s'}$ and $j_2$ be the continuous inclusion $V^{s+2} \hookrightarrow V^{s'+2}$.  Then $j_2\circ Y_s$ is  continuous as a $V^{s'}$-valued function on $[0, \infty)$, differentiable as a $V^{s'}$-valued function on $(0, \infty)$, takes values in $V^{s'+2}$, and
$$\tfrac{d}{dt} \,\iota'\circ ( j_2 \circ Y_s) \;\; = \;\; j_0 \circ \tfrac{d}{dt}(\iota \circ Y_s) \;\; = \;\; j_0 \circ (M_s \circ Y_s) \;\; = \;\; M_{s'} \circ (j_2 \circ Y_s) \, , $$
i.e. $j_2 \circ Y_s$ satisfies (\ref{sprime-IVP}). By Proposition \ref{uniqueness-abstract}, $j_2 \circ Y_s = Y_{s'}$, and thus for all $t$, $Y_s(t) = j_2\circ Y_s(t) = Y_{s'}(t)$.

Note that this also shows that for any $s \leq -\ell-2$, there is no other function besides $\widetilde{U}$ that is continuous as a $V^s$-valued function on $[0, \infty)$, is differentiable as a $V^s$-valued function on $(0, \infty)$, takes values in $V^{s+2}$, and solves (\ref{IVP}). \end{proof}
\end{prop}

\begin{remark*}  Our two constructions of the spectral coefficient function, explicitly in  \textsection\textsection\ref{spectral-constr} and via semigroup theory in the proof just above, employ contrasting ways of thinking about multiplication by a function in weighted $L^2$-spaces.  In the first perspective, the weighted $L^2$-spaces $V^s$, $s \in \Z$, form a nested family of Hilbert spaces, with dense continuous inclusions, and the multiplication maps $(1-\mu)$ are continuous linear maps $V^{s+2} \to V^s$.    In the second perspective, each $V^{s+2}$ is considered as the dense subspace of $V^s$, and each multiplication map is a densely defined unbounded operator, which is the infinitesimal generator for a strongly continuous semigroup.  To show that the two constructions yield the same function in the proof of the preceding proposition, we needed to consider not just one multiplication operator $V^{s+2} \to V^s$ but the \emph{family} of multiplication operators on the nested family of weighted $L^2$-spaces.  
\end{remark*}

Recall from  \textsection \ref{spectral-constr} that we construct an automorphic heat kernel, $U$, from $\widetilde{U}$ by spectral synthesis in Theorem \ref{constr-afc-heat-ker}.  We can now show that $U$ is the unique automorphic heat kernel.

\begin{theorem}\label{uniqeness-afc-heat-ker}   The automorphic heat kernel constructed in Theorem \ref{constr-afc-heat-ker} is the \emph{unique} automorphic heat kernel, as defined in  \textsection\ref{more-precise}.  Moreover it is strongly continuously differentiable as a $H^{-\ell-2}$-valued function on $[0, \infty)$.
\begin{proof} This follows immediately from Lemma \ref{translation} and Proposition \ref{uniqueness-coef}.\end{proof}
\end{theorem}

\section{Smoothness of the Automorphic Heat Kernel}\label{smoothness}

Our last result is the smoothness of the automorphic heat kernel as a function on $X$, for fixed $t > 0$.  We again work on the spectral side, showing that for $t>0$, $\widetilde{U}(t)$ lies in $V^s$, for \emph{every} index $s$.  We will use this to show that, for $t >0$, the automorphic heat kernel takes values in $H^\infty$ and thus in $C^\infty$, by global automorphic Sobolev embedding.

\begin{prop} \label{V-infinity} For $t>0$, for all $s$, $\widetilde{U}(t) \in V^s$, i.e. $\widetilde{U}(t) \in V^\infty$.
\begin{proof} We must show  $(1-\lambda_\xi)^{\nicefrac{s}{2}} \overline{\Phi}_\xi(x_0) e^{\lambda_\xi t} \,  \in \, L^2(\Xi)$, for all $s$.  Since $\mathcal{F} \delta \in V^{-\ell}$, we know $(1-\lambda_\xi)^{\nicefrac{-\ell}{2}} \, \overline{\Phi}_\xi(x_0) \, \in \, L^2(\Xi)$.  Thus it suffices to show that, for $m = s + \ell$,
$$(1-\lambda_\xi)^{\nicefrac{m}{2}}\,  e^{\lambda_\xi t} \;\; \in \;\; L^\infty(\Xi) \;.$$
Note that this is a continuous function of $\xi$ since the map $\xi \mapsto \lambda_\xi$ is continuous.  To see that this function is bounded, consider $\Xi$ as an ascending union of the compact sets $\Xi_T = \{\xi \in \Xi \, : \, |\lambda_\xi| \leq T\} \;.$  Let $M_T$ be the maximum on $\Xi_T$, i.e. 
$$M_T \;\ = \;\; \max_{\xi \, \in \, \Xi_T} \;  (1-\lambda_\xi)^{\nicefrac{m}{2}}\,  e^{\lambda_\xi t}\;\;= \;\;\max_{\xi \, \in \, \Xi_T} \;  (1+|\lambda_\xi|)^{\nicefrac{m}{2}} \, e^{-t |\lambda_\xi|}\;.$$
Then $\displaystyle{\lim_{T \to \infty}} M_T \, < \, \infty$, since
$$| M_{T+1} \, - \, M_T|  \; \leq \; \max_{\xi \;\in \;\Xi_{T+1}-\Xi_T} \;  (1+|\lambda_\xi|)^{\nicefrac{m}{2}} \, e^{-t |\lambda_\xi|} \; < \; (2+T)^{\nicefrac{m}{2}} \, e^{-tT}  \xrightarrow{T \to \infty}  0 \;,$$
and $(1+|\lambda_\xi|)^{\nicefrac{m}{2}} \, e^{-t |\lambda_\xi|}$ is bounded on $\Xi$. \end{proof}
\end{prop}

\begin{theorem} For $t>0$, the automorphic heat kernel lies in $C^\infty(X)$, and its automorphic spectral expansion
$$U(t) \;\; = \;\; \int_\Xi \overline{\Phi}_\xi(x_0) \cdot e^{\lambda_\xi\,  t} \cdot \Phi_\xi \, d\xi $$
converges in the $C^\infty(X)$-topology.  
\begin{proof}  This follows from Proposition \ref{V-infinity}, the fact that $\mathcal{F}^{-1}$ maps $V^{s}$ to $H^{s}$ for all $s$, and the global automorphic Sobolev embedding theorem. \end{proof}
\end{theorem}


\appendix
\section{Limits and Derivatives of Vector-valued Functions}\label{lim-deriv-V-vald}

We collect results about limits and differentiation of vector-valued and operator-valued functions which are needed in the proofs above or in \textsection\textsection\ref{acp-uniqueness}, below.

\begin{lemma}\label{comp-w-diffble}  For topological vector spaces $V$ and $W$, a continuous linear map $T:V\to W$, and a $V$-valued function $f$ of a real variable, if $f$ is weakly (or strongly) differentiable, then $T \circ f$ is weakly (resp. strongly) differentiable, as a $W$-valued function, and its weak (resp. strong) derivative satisfies $(T \circ f)' = T \circ f'$.  Moreover if $f$ is weakly (or strongly) continuously differentiable, then $T \circ f$ is also.
\begin{proof} The proof is a straightforward application of the definitions. \end{proof}
\end{lemma}

Let $V$ be a Banach space and $\textbf{B}(V)$ the space of bounded linear operators on $V$.

\begin{lemma}\label{limits-op-v}  Let $\mathcal{I}$ be a real interval and $\Theta: \mathcal{I} \to  \textbf{B}(V)$ be a strongly continuous operator-valued function that is bounded on compact intervals: for any $[\alpha, \beta] \subset \mathcal{I}$, there is a number $C$ such that $\lVert \Theta(t) \rVert_{\mathrm{op}} \leq C$ for all $t \in [\alpha, \beta]$.   Let $a \in \mathcal{I}$, and let  $v(t)$ be a $V$-valued function on $\mathcal{I}$ such that the limit of $v(t)$ as $t\to a$ exists.  Then 
$$\lim_{t \to a} \Theta(t) v(t) \;\; = \;\; \Theta(a) \big(\lim_{t \to a} v(t) \big) \;,$$
where both limits are in the topology of $V$.
\begin{proof}  This is a standard two epsilon proof. \end{proof}

\end{lemma}

\begin{remark*} For a \emph{uniformly continuous} operator-valued function $\Theta(t)$, the hypothesis that $\Theta$ is bounded on compacts is unnecessary and the result is immediate.  However, we apply the lemma to translations of strongly continuous semigroups, which are not necessarily uniformly continuous.  In fact, a  semigroup is uniformly continuous if and only if it has a bounded (everywhere-defined) infinitesimal generator.  (See e.g. \cite{dunford-schwartz1958}.)\end{remark*}

\begin{lemma}\label{limits-v-c}  Let $v(t)$ be a $V$-valued function such that the limit of $v(t)$ as $t \to a$ exists and $\alpha(t)$ a scalar-valued function such that the limit of $\alpha(t)$ as $t \to a$ exists.  Then the limit of $\alpha(t) v(t)$ as $t \to a$ exists in $V$ and 
$$\lim_{t \to a} \big(\alpha(t) \, v(t) \big) \;\; = \;\; \big(\lim_{t \to a} \alpha(t) \big) \cdot \big(\lim_{t \to a} v(t) \big) \;.$$
\begin{proof}  This is another standard two epsilon proof. \end{proof}
\end{lemma}

\begin{lemma}[Generalized Chain Rule] \label{chain-rule} Let $v$ be a strongly differentiable $V$-valued function on an open interval $\mathcal{I}$ and $g(t)$ a differentiable monotonic scalar-valued function, taking values in $\mathcal{I}$. Then $t \mapsto  v(g(t))$ is a strongly differentiable $V$-valued function on the domain of $g$, and its derivative is 
$$\tfrac{d}{dt} v(g(t)) \;\;= \;\; v'(g(t)) \cdot g'(t) \;.$$
\begin{proof} Lemma \ref{limits-v-c} allows us to prove this as for scalar-valued functions. \end{proof}

\end{lemma}

We recall the \emph{weak-to-strong differentiability} theorem, which is a consequence of the weak-to-strong boundedness principle.  These results hold very generally, for $V$ any quasi-complete locally convex vector space.  See \cite{garrett2018}, Chapter 15.

\begin{theorem*}[15.1.1 in \cite{garrett2018}] For $k \geq 1$, a weakly $C^k$ $V$-valued function on a real interval is  is strongly $C^{k-1}$.\end{theorem*}

\section{Vector-valued Integrals (Gelfand-Pettis)}\label{gelfand-pettis}

We used the vector-valued (weak) integrals of Gelfand \cite{gelfand36} and Pettis \cite{pettis38} in the proof of Proposition \ref{spectral-diff-eqtn}.  Here we describe the Gelfand-Pettis integral and summarize the key results; see \cite{garrett2018}, 14.1.

Let $X$ be a \emph{compact} Hausdorff topological space with \emph{finite} positive regular Borel measure $\mu$.  Let $V$ be a locally convex, quasi-complete topological vector space.  The \emph{Gelfand-Pettis integral} is a vector-valued integral $C^0(X,V) \to V$ denoted:
$$f \;\; \mapsto \;\; I_f \; = \; \int_X f \, d\mu \; , $$
characterized by the property that, for all $\lambda \in V^{\ast}$, 
$$ \lambda\left(\int_X f \, d\mu\right) \; = \; \int_X \lambda \circ f \; d\mu \;,$$
where this latter integral is the usual scalar-valued Lebesgue integral.  

\begin{remark*} Alternately we may allow $X$ to be non-compact, as long as it is locally compact and of finite measure; in this case we require $f$ to be compactly supported. \end{remark*}

\begin{remark*} (i) The precise hypothesis on $V$ that is necessary for the existence of the Gelfand-Pettis integral is that the closure of the convex hull of a compact set is compact.  Quasicompleteness ensures this.  (ii) Hilbert, Banach, Frechet, LF spaces, and their weak duals are locally convex, quasi-complete topological vector spaces; see \cite{garrett2018}, Chapter 13.
\end{remark*}

\begin{theorem*}\label{gp_ints}
\textrm{(i)} The Gelfand-Pettis integral \emph{exists}, is \emph{unique}, and satisfies the following \emph{estimate}:
$$ I_f \; \in \; \mu(X) \cdot \big(\text{closure of convex hull of } f(X) \big) \; .$$
\textrm{(ii)} Any continuous linear operator $T$ from $V$ to another locally convex, quasi-complete topological vector space $W$ commutes with the Gelfand-Pettis integral: 
$$T \left(\int_X f \, d\mu\right) \;\; = \;\; \int_X T\circ f \, d\mu \;.$$
\end{theorem*}

\section{Unbounded Operators on Hilbert Spaces}\label{unbdd-ops}

In this appendix, we give a concise recounting of the definitions and results from the theory of unbounded operators needed for the discussion of multiplication operators on weighted $L^2$-spaces in  \textsection\ref{mult-ops}, above.  See Chapter 9 in Garrett's book \cite{garrett2018} or Chapter 12 in Grubb's book \cite{grubb2009}, for more extensive treatments.

Let $V$ be a Hilbert space and $D$ a subspace.  A linear map $T:D \to V$ is called an \emph{unbounded operator on $V$}.  Specifying the domain $D$ is an essential part of defining an unbounded operator.  Because of this and for brevity, we sometimes introduce an unbounded operator as a pair $T, D$.  The domain of an unbounded operator $T$ is also denoted $\mathrm{Dom}(T)$.

An unbounded operator on $V$ is \emph{closed} (or \emph{graph-closed}) if its graph is closed in $V \oplus V$.  For everywhere-defined linear operators the notion of closedness coincides with that of continuity: a continuous linear operator has a closed graph, and, by the Closed Graph Theorem, an everywhere-defined linear operator with a closed graph is continuous.   In contrast, unbounded operators with closed graphs are \emph{not} necessarily continuous.   

An unbounded operator $T_2$ is an \emph{extension of} an unbounded operator $T_1$ if it has a larger domain and it agrees with $T_1$ on the domain of $T_1$, i.e.
$$\mathrm{Dom}(T_2) \, \supset \, \mathrm{Dom}(T_1) \;\;\;\;\;\;\;\; \text{and}\;\;\;\;\;\;\;\;T_2 |_{\mathrm{Dom}(T_1)} \; = \; T_1 \;.$$
The expression $T_2 \supset T_1$ denotes that $T_2$ is an extension of $T_1$.  An unbounded operator $T', D'$ is a \emph{subadjoint} to $T,D$, when 
$$\langle Tv, w \rangle = \langle v, T'w\rangle \;\;\;\;\;\;\;\; \text{ for all } v \in D, \; w \in D' \;.$$  The \emph{adjoint} $T^\ast$ of $T$ is defined as the unique maximal element among all subadjoints.  It is only guaranteed to exist if $T$ is densely defined; even so existence and uniqueness of the adjoint require proof.

\begin{prop*}[9.1.1 in \cite{garrett2018}]  Let $T, D$ be an unbounded, densely defined operator on a Hilbert space $V$.  \begin{enumerate}
\item[(i)] There is a unique maximal $T^\ast, D^\ast$ among all subadjoints to $T,D$.
\item[(ii)] $T^\ast$ is closed, in the sense that its graph is closed in $V \oplus V$.
\item[(iii)] $T^\ast$ is characterized by its graph:
$$\mathrm{Graph}(T^\ast) \; = \; \big(U(\mathrm{Graph}(T))\big)^\perp\;,$$
where $U$ is the Hilbert space isomorphism $V \oplus V \to V \oplus V$ given by $U(v \oplus w) \; = \; -w \oplus v$.
\end{enumerate}

\end{prop*}

\begin{remark}\label{dom-adjt} The proof of this result shows that the domain of the adjoint $T^\ast$ of $T$ consists precisely of those vectors $w \in V$ for which there exists $w' \in V$ such that $\langle Tv, w \rangle = \langle v, w'\rangle$ for all $v \in \mathrm{Dom}(T)$.  \end{remark}

For bounded operators, being symmetric and being self-adjoint are equivalent; however in the case of unbounded operators the notions are distinct.  For an unbounded operator $T$ to be \emph{symmetric} means:
$$\langle Tv, w \rangle \; = \; \langle v, Tw \rangle \;\;\;\;\;\;\; \text{for all } v,w \in \mathrm{Dom}(T)  \;.$$
In the case where $T$ is densely defined, so $T^\ast$ is guaranteed to exist, an equivalent criterion is that $T \subset T^\ast$.  A densely defined unbounded operator is \emph{self-adjoint} when $T = T^\ast$, i.e. when (i)  the domain of its adjoint, which is maximal among domains of subadjoints, is \emph{precisely} the domain of $T$ and (ii) $T$ and its adjoint agree on the domain of $T$.  Clearly, a self-adjoint operator is symmetric.  Note that it is also closed, since adjoint operators are closed.  

\begin{remark}\label{dd-sym-self-adjt} A densely defined symmetric operator $T$ is self-adjoint if its domain contains the domain of its adjoint, since in this case we have both $T \subset T^\ast$ (by symmetry) and $T \supset T^\ast$ and thus $T = T^\ast$.\end{remark}

A symmetric operator $T$ is \emph{lower semi-bounded} if there is a real constant $c$ such that  $\langle T v, v \rangle \; \geq \; c \, \langle v, v \rangle$ for all $v \in \mathrm{Dom}(T)$; \emph{positivity} is a special case of lower semi-boundedness, with $c=0$.  A symmetric operator is \emph{upper semi-bounded} if there is a real constant $C$ such that $\langle T v, v \rangle \; \leq \; C \, \langle v, v \rangle$  for all $v \in \mathrm{Dom}(T) $; \emph{negativity} is a special case of upper semi-boundedness with $C =0$.

The \emph{resolvent set} of an unbounded operator is the set of $\lambda \in \C$ for which the operator $T-\lambda I$ is a bijection of $\mathrm{Dom}(T)$ onto $V$ with bounded inverse $(T - \lambda I)^{-1}$.

\section{Semigroup Theory and the Abstract Cauchy Problem}\label{semigroup-theory}

We summarize the results of semigroup theory that were used in  \textsection \ref{uniqueness-section} to prove the uniqueness of the automorphic heat kernel, as defined precisely in  \textsection\textsection\ref{more-precise}.

\subsection{Semigroups of linear operators}

We mostly follow Grubb's account in Sections 14.2 and 14.3 of \cite{grubb2009}, on contraction semigroups in Banach spaces and Hilbert spaces, respectively.  Classic references for the theory of semigroups are the books of Hille and Phillips \cite{hille-phillips1957} and Dunford and Schwartz \cite{dunford-schwartz1958}.  For more modern treatments see, e.g., \cite{pazy1983} or \cite{goldstein1985}.

Let $V$ be a Banach space and $\mathbf{B}(V)$ the space of bounded linear operators on $V$.  For a function $G:[0, \infty) \to \mathbf{B}(V)$ to be a \emph{semigroup} means: 
$$G(0) \; =\; I \hskip 1cm \textrm{and} \hskip 1cm G(s+t) \;\; = \;\; G(s) \cdot G(t) \hskip .5cm \text{for all } s, t \geq 0 \;.$$
For a semigroup $G(t)$ to be \emph{strongly continuous} means: for each $v \in V$, the map $t \mapsto G(t) v$ is continuous as a $V$-valued function on $[0, \infty)$, i.e. $G$ is strongly continuous as an operator-valued function.  For $G(t)$ to be a \emph{contraction} semigroup means:
$$\lVert G(t) \lVert_{\mathrm{op}}\;\; \leq \;\; 1 \hskip .5cm \text{for all } t \geq 0 \;.$$
The \emph{infinitesimal generator} of $G(t)$ is the unbounded operator $A$ on $V$ defined by 
$$A v  \;\; = \;\; \lim_{h \to 0^+} \; \tfrac{1}{h}\big(G(h) - I\big)v \; ,$$
the domain consisting precisely of those $v$ for which the limit exists.  

\begin{remark*} A \emph{contraction} semigroup $G(t)$ is strongly continuous if it is strongly continuous at $t=0$.  See Lemma 14.1 in \cite{grubb2009}.  For this reason, Grubb \emph{defines} a strongly continuous contraction semigroup to be a contraction semigroup with this additional property. \end{remark*}

\begin{remark*} In fact, \emph{any} semigroup $G(t)$ that is strongly continuous at $t=0$ is strongly continuous on $[0, \infty)$.  For this reason, some authors only require strong continuity at zero in the definition of strongly continuous semigroup. \end{remark*}

While not every strongly continuous semigroup $G(t)$ is a contraction semigroup, we do have a bound for the operator norm of $G(t)$ in terms of $t$, as follows.

\begin{prop*}[Theorem 2.2 in \cite{pazy1983}] Given a strongly continuous semigroup $G(t)$, there are nonnegative constants $M$ and $\omega$ such that 
$$\lVert G(t) \lVert_{\mathrm{op}} \;\; \leq \;\; Me^{\omega t} \hskip .5cm \text{for all } t \geq 0 \;.$$
\end{prop*}

\begin{remark}\label{bdd-finite-intvls} As a consequence, strongly continuous semigroups are bounded on compact intervals: if $G(t)$ is a strongly continuous semigroup, then for any $\tau \geq 0$ there is a constant $C$ such that $\lVert G(t) \rVert_{\mathrm{op}} \leq C$ for all $t \in [0, \tau]$.
\end{remark}

\begin{theorem*}[Theorem 14.2, \cite{grubb2009}] Let $G(t)$ be a strongly continuous contraction semigroup with infinitesimal generator $A$, having domain $D$. For $v_0 \in D$,  the vector-valued function $t \mapsto G(t) v_0$ is strongly differentiable, takes values in $D$, and satisfies
$$\tfrac{d}{dt}\big(G(t) v_0\big) \;\; = \;\; G(t) \, Av_0 \;\; = \;\; A \big(G(t) v_0\big) \hskip .5cm \textrm{for all } t \geq 0 \, .$$
\end{theorem*}

\begin{remark}\label{Grubb-14.2-extd} This theorem holds for \emph{any} strongly continuous semigroup, not only those consisting of contractions.  See, e.g. Theorem 10.3.3 in \cite{hille-phillips1957}, although the proof provided there omits some details.  Grubb's careful proof of her Theorem 14.2 can be modified to prove the analogous result for any strongly continuous semigroup: the property of boundedness on compact intervals (see Remark \ref{bdd-finite-intvls}) can be used instead of the contraction property.\end{remark}

\begin{remark}\label{Grubb-14.2-C1}  The conclusion of the theorem can be strengthened also: $t \mapsto G(t) v_0$ is $C^1$, not merely strongly differentiable, since $Av_0 \in V$ implies $t \mapsto G(t) A v_0$ is continuous, by the strong continuity of $G$. \end{remark}

The converse problem, under what conditions a given unbounded operator is the infinitesimal generator of a strongly continuous semigroup, was solved independently by Hille and Yosida in 1948 \cite{hille1948, yosida1948}.  We use the following consequence in the proof of Proposition \ref{uniqueness-coef}.

\begin{theorem*}[Corollary 14.11, \cite{grubb2009}]  An unbounded operator $A$ on a Hilbert space $H$ is the infinitesimal generator of a strongly continuous contraction semigroup if and only if $A$ is densely defined, closed, upper semi-bounded with upper bound $C \leq 0$ and has the positive real numbers contained in its resolvent set. \end{theorem*}

\subsection{The abstract Cauchy problem: uniqueness of solutions} \label{acp-uniqueness}

The abstract Cauchy problem (ACP) has been formulated in various ways (see, e.g.  \textsection23.7 in \cite{hille-phillips1957}, \cite{phillips1954}, \cite{beals1972},  \textsection14.4 in \cite{grubb2009},  \textsection3.1 in \cite{arendt-etal2001}), but roughly it is to find a vector-valued solution $Y(t)$ of an initial value problem of the form
$$\tfrac{d}{dt} Y(t) \;\; = \;\; A Y(t) \; ; \hskip .5cm Y(0) = v_0 \;,$$
where $A$ is an unbounded operator on a Banach space $V$ and $v_0 \in V$.  Different formulations of the ACP impose various additional conditions.

As discussed in the previous subsection, semigroup theory provides a solution in the case that $A$ is the infinitesimal generator of a strongly continuous semigroup and $v_0$ is in the domain of $A$.  This solution is continuously differentiable and takes values in the domain of $A$.  In Proposition \ref{uniqueness-abstract} below, we formulate a version of the abstract Cauchy problem and show that it has a unique solution.  

\begin{remark*} Although there are many theorems on the uniqueness of solutions to the ACP in the literature (see e.g. \cite{phillips1954}, \cite{beals1972}, Theorem 4.1.3 in \cite{pazy1983}, Theorem II.1.2 in \cite{goldstein1985}), some use formulations of the ACP which are not suitable for our application and others omit details in their proof.  Thus, for the sake of completeness, we include a careful formulation of the ACP and proof of the uniqueness of its solution. \end{remark*}

We will need the following lemma, which can perhaps be viewed as a generalization of the Leibntiz rule.

\begin{lemma}\label{prod-rule}  Let $G(t)$ be a strongly continuous semigroup with infinitesimal generator $A$ having domain $D$.  For any $v(t)$ strongly differentiable on an open interval $\mathcal{I} \subset (0, \infty)$ and taking values in $D$, the function $t \mapsto G(t) v(t)$ is strongly differentiable on $\mathcal{I}$, takes values in $D$, and its derivative is:
$$\tfrac{d}{dt} \big( t \mapsto G(t) \, v(t) \big) \;\; = \;\; t \;  \mapsto \; G(t) v'(t) \, + \, G(t) Av(t) \;.$$
\begin{proof} We will consider the left and right difference quotients separately.  We will use Lemma \ref{limits-op-v} to evaluate limits involving both operator-valued and vector-valued functions.  Note that, since $G$ is a strongly continuous semigroup, it is strongly continuous as an operator-valued function and also bounded on compact intervals, by Remark \ref{bdd-finite-intvls}.  So then, for any real number $h$, the operator-valued function $h \mapsto G(t+h)$, although not a semigroup, is also strongly continuous and bounded on compact intervals, i.e. it satisfies the hypotheses of Lemma \ref{limits-op-v}.

First consider the difference quotient for $h >0$.  Rewrite it as in the proof of the Leibnitz rule for  scalar-valued functions, and then make use of one of the semigroup properties of $G$, namely that for  \, $t, h \geq 0$, \, $G(t+h) \, = \, G(t) \,  G(h)$, to obtain:
\begin{eqnarray*}
\lefteqn{\tfrac{1}{h} \big(G(t+h)v(t+h) - G(t) v(t) \big)}\\
& = &  G(t+h)\bigg[ \;\tfrac{1}{h} \;  \big(v(t+h)-v(t)\big) \bigg] \;\; +\;\; G(t) \; \bigg[\; \tfrac{1}{h} \; \big(G(h)-I\big)\; \bigg]\;v(t) \;.
\end{eqnarray*}
Applying the limit as $h \to 0^+$ to the first term, since $v$ is strongly differentiable, Lemma \ref{limits-op-v} implies
$$\lim_{h \to 0^+} \; G(t+h)\bigg[ \;\tfrac{1}{h} \;  \big(v(t+h)-v(t)\big) \bigg]  \;\; = \;\; G(t) v'(t) \;.$$
Applying the limit to the second term, we use the fact that $G(t)$ is continuous on $V$, the fact that $v(t) \in D$,  and the definition of the infinitesimal generator to conclude:
$$\lim_{h \to 0^+} \;  G(t) \; \bigg[\; \tfrac{1}{h} \; \big(G(h)-I\big)\; \bigg]\;v(t) \;\; = \;\; G(t) A v(t) \;.$$
Note that since $G(t)$ and $G(h)$ commute, this also shows that the limit, 
$$\lim_{h \to 0^+} \; \bigg[\; \tfrac{1}{h} \; \big(G(h)-I\big)\; \bigg]\;  G(t) \, v(t) \;, $$
exists, i.e. that $G(t) \, v(t)$ is in the domain of $A$ for all $t$.  Thus $t \mapsto G(t) v(t)$ takes values in $D$, as claimed.  Moreover, we have shown
$$\lim_{h \to 0^+} \tfrac{1}{h} \big(G(t+h)v(t+h) - G(t) v(t) \big) \;\; = \;\; G(t) v'(t)  \;+ \;G(t) A v(t) \;.$$
Next we consider the difference quotient for $h<0$,
\begin{eqnarray*}
\lefteqn{\tfrac{1}{h} \big(G(t+h)v(t+h) - G(t) v(t) \big)}\\
& = &  G(t+h)\bigg[ \;\tfrac{1}{h} \;  \big(v(t+h)-v(t)\big) \bigg] \;\; +\;\; \tfrac{1}{h} \; \bigg[ \big(G(t+h)-G(t)\big)v(t) \; \bigg] \;.
\end{eqnarray*}
The limit of the first term as \, $h \to 0^-$ \, is \, $G(t) v'(t)$.  Rewriting the second term, with $\eta = -h$, and using the semigroup property, the second term becomes:
$$ \tfrac{1}{\eta} \; \bigg[ \big(G(t-\eta+\eta)-G(t-\eta)\big)v(t) \; \bigg]\;\;= \;\;  G(t-\eta) \; \bigg[ \; \tfrac{1}{\eta}  \big(G(\eta)-I\big) \; \bigg] \; v(t) \;.$$
Taking the limit as $\eta \to 0^+$,  we obtain 
$$ \lim_{\eta \to 0^+} \; G(t-\eta) \; \bigg[ \; \tfrac{1}{\eta}  \big(G(\eta)-I\big) \; \bigg] \; v(t)\;\; = \;\; G(t) \, A v(t) \; ,$$
using the strong continuity of $G$, Lemma \ref{limits-op-v}, the fact that $v(t) \in D$, and the definition of the infinitesimal generator.  Thus we have shown
$$\lim_{h \to 0-} \; \tfrac{1}{h} \big(G(t+h)v(t+h) - G(t) v(t) \big)\;\; = \;\; G(t) v'(t) \;+ \; G(t) \, A v(t) \;.$$
Therefore $ t \mapsto G(t) v(t)$ is strongly differentiable,  with derivative \, $t \; \mapsto \; G(t) v'(t) \,+ \, G(t) \, A v(t) $, as claimed, and takes values in $D$.  \end{proof}
\end{lemma}

\begin{prop} \label{uniqueness-abstract} Let $G(t)$ be a strongly continuous semigroup in a Banach space $V$, let $A$ be the infinitesimal generator for $G(t)$ with domain $D$, and let $v_0$ be an element of $D$.  Then there is a unique function $[0, \infty) \to V$ that (i) is strongly continuous on $[0, \infty)$, (ii) is strongly differentiable on $(0,\infty)$, (iii) takes values in $D$, and (iv) solves the initial value problem, 
$$\tfrac{d}{dt} Y(t) \;\; = \;\; A \; Y(t) \; ; \hskip .5cm Y(0) \; = \;\; v_0\;.$$
Moreover, the solution is strongly continuously differentiable on $[0, \infty)$.

\begin{proof}
Let $Y_1(t)$ be the known solution $t \mapsto G(t) v_0$, which is in fact $C^1$, and $Y_2(t)$ another strongly differentiable solution taking values in $D$.  Fix $\tau > 0$.  We will show that $Y_1(\tau) = Y_2(\tau)$.  Let $w(t) = G(\tau - t) Y_2(t)$, so that
\begin{eqnarray*}
w(0) & = &  G(\tau) Y_2(0) \; = \; G(\tau) v_0 \; = \; Y_1(\tau) \;,  \textrm{ and} \\
w(\tau) & = &  G(0) Y_2(\tau) \;= \; I \, Y_2(\tau) \; = \; Y_2(\tau) \;.
\end{eqnarray*}
We claim that $w$ is strongly continuous on $[0,\tau]$.  Indeed, since $G$ is a strongly continuous semigroup, it is strongly continuous as an operator-valued function and also bounded on compact intervals, by Remark \ref{bdd-finite-intvls}.  Thus the operator-valued function $t \mapsto G(\tau-t)$, although not a semigroup, is also strongly continuous and bounded on compact intervals.  Applying Lemma \ref{limits-op-v} proves that $w$ is strongly continuous on $[0, \tau]$.

We show that $w$ is strongly differentiable with identically zero derivative on $(0, \tau)$.  Let $v(t) = Y_2(\tau -t)$.  Then $v$ takes values in $D$ and, by Lemma \ref{chain-rule}, is strongly differentiable on $(0, \tau)$ with derivative $v'(t) = -Y_2'(\tau-t)$.   By Lemma \ref{prod-rule}, the function $t \mapsto G(t) v(t)$ is strongly differentiable on $(0, \tau)$, takes values in $D$, and has derivative $t \mapsto G(t) v'(t) + G(t) Av(t)$.  Viewing $w(t)$ as the composition of $t \mapsto G(t) v(t)$ with $t \mapsto (\tau-t)$ and applying Lemma \ref{chain-rule} once again, we see that $w(t)$ is strongly differentiable on $(0, \tau)$, takes values in $D$ and has derivative 
$$w'(t) \;\; = \;\; - \big(G(\tau-t)v'(\tau-t) \, + \, G(\tau-t)Av(\tau-t)\big) \;\; = \;\; G(\tau-t)Y_2(t) \, - \, G(\tau-t)AY_2(t) \;.$$
Since $Y_2$ satisfies the differential equation, $Y_2'(t) = A \, Y_2(t)$, and thus $w'(t)$ is identically zero on $(0, \tau)$.

Take any $v^\ast$ in $V^\ast$.  Then the scalar-valued function $t \mapsto \langle w(t), v^\ast \rangle$ is continuous on $[0,\tau]$ and has identically zero derivative on $(0, \tau)$, so is constant on $[0, \tau]$.  In particular,
$$\langle w(0), v^\ast \rangle \;\; = \;\; \langle w(\tau), v^\ast \rangle \hskip .5cm \textrm{i.e.} \hskip .5cm \langle Y_1(\tau) , v^\ast \rangle \;\; = \;\; \langle Y_2(\tau), v^\ast\rangle \;.$$
Since this holds for all $v^\ast \in V^\ast$, the Hahn-Banach theorem implies that $Y_1(\tau) = Y_2(\tau)$.  Further, since $\tau$ was an arbitrary positive number, we see that $Y_1$ and $Y_2$ agree on $(0, \infty)$.  Since $Y_1(0) = v_0 = Y_2(0)$, we have $Y_1 = Y_2$ on $[0, \infty)$. \end{proof}
\end{prop}

\bibliography{afc-heat-ker-biblio}{}
\bibliographystyle{plain}

\end{document}